\documentclass{article}

\usepackage{amsmath}
\usepackage[latin1]{inputenc}
\usepackage[T1]{fontenc}
\usepackage[english]{babel}
\usepackage{enumerate}
\usepackage{amssymb}
\usepackage{amsthm}
\usepackage{cleveref}
\usepackage{lipsum}

\newtheorem{Lemma}{Lemma}[section]
\newtheorem{Theorem}[Lemma]{Theorem}
\newtheorem{Corollary}[Lemma]{Corollary}
\newtheorem{Remark}[Lemma]{Remark}
\newtheorem{Example}[Lemma]{Example}

\title{Spectra of typical Hilbert space operators}
\author{Marcel Scherer}
\begin{document}
\maketitle
\begin{abstract}Let $\mathcal{B}(H)$ denote the algebra of all bounded linear operators on a separable Hilbert space, equipped with the norm topology. A property is called typical if the set of operators fulfilling the property is co-meager. We show that having non-empty continuous spectrum is a typical property and that the set of operators with empty continuous spectrum is dense in $\mathcal{B}(H)$. In addition, we show that the set of operators with empty point spectrum is nowhere dense. Moreover, we characterize the closure of the set of operators whose spectrum and point spectrum coincide.
\end{abstract}
\section{Introduction}
\let\thefootnote\relax\footnotetext{
2020 \textit{Mathematics Subject Classification.} 47A11.\\
\ \textit{Key words and phrases.} $G_\delta$-sets, co-meager, point spectrum, continuous spectrum, (Semi)- Fredholm operators, polar decomposition, Banach-Mazur game, Similarity Theorem.\\
 The author was partially supported by the Emmy Noether Program of the German Research Foundation (DFG Grant 466012782).}

$G_\delta$-sets appear in many different contexts, for example, as the continuity points of a real valued function or as the extreme points of any metrizable convex compact set in a locally convex topological vector space. A subspace of a Baire space $X$ is called \textit{co-meager} if it contains a dense $G_\delta$-set and a property $\Phi$  is called typical if $\{x\in X,\ x\ \textup{has property}\ \Phi\}$ is co-meager. Typical properties of bounded linear operators on certain Banach spaces were recently studied in \cite{GMM1}, \cite{GMM2} and \cite{GM}. Furthermore, T. Eisner and T. M\'{a}trai \cite{TE} studied typical properties of operators on Hilbert spaces with respect to the weak, strong, strong-star and the norm topology. In this paper we will only study typical properties with respect to the operator norm topology and answer the two questions given in \cite[Problem 8.4]{TE}.\\
We start with the basic definitions of the point and continuous spectrum in \Cref{Section 2}, followed by some well-known results about semi-Fredholm operators and helpful lemmas, which we need for the proofs of the main theorems.\\
In \Cref{Section 3}, we look at the first question in \cite[Problem 8.4]{TE}: Is the set of operators with non-empty point spectrum co-meager in $\mathcal{B}(H)$ with respect to the norm topology? We answer this question in the affirmative by showing the even better result:
\begin{Theorem}
Let $H$ be a separable Hilbert space. Then the set 
  \begin{equation*}
    \{T\in\mathcal{B}(H);\ \sigma_p(T)=\emptyset\}
  \end{equation*}
is nowhere dense with respect to the norm topology.
\end{Theorem}
In \Cref{Section 4} we move on to the second question of \cite[Problem 8.4]{TE} which asks if the set of operators with empty continuous spectrum is co-meager. Unlike the previous question, the answer is $\grqq$no$\grqq$ and we will prove:
\begin{Theorem}
Let $H$ be a separable infinite dimensional Hilbert space. Then the set
\begin{equation*}
\{T\in\mathcal{B}(H);\ \sigma_c(T)\neq\emptyset\}
\end{equation*}
is co-meager with respect to the norm topology.
\end{Theorem}
The proof of this theorem uses the Banach-Mazur game and will be somewhat technical.\\
In the last section we will answer a question that arose during the study of operators with empty continuous spectrum: What is the closure of the set of operators with empty continuous spectrum?
Knowing that this set is meager suggests that its closure is not very large, but using the Similarity Orbit Theorem of Apostol, Fialkow, Herrero and Voiculescu \cite[Theorem 9.2]{DH2} we will show that it is already dense.
\begin{Theorem}
Let $H$ be a separable Hilbert space. Then the set 
  \begin{equation*}
    \{T\in\mathcal{B}(H);\ \sigma_c(T)=\emptyset\}
  \end{equation*}
is dense in $\mathcal{B}(H)$ with respect to the norm topology.
\end{Theorem}

\section*{Acknowledgement}
I wish to thank Michael Hartz for many helpful discussions and John McCarthy for pointing me towards the Similarity Orbit Theorem. I would also like to thank the referee for the thorough reading and for providing historical context and additional references regarding the proof of \Cref{MainTheorem3}.

\section{Definitions}
\label{Section 2}
Throughout the paper $H$ will be a separable Hilbert space and $\mathcal{B}(H)$ the set of all linear, bounded operators from $H$ to $H$. For $T\in\mathcal{B}(H)$, the \textit{spectrum $\sigma(T)$ of T} are all points $\lambda\in\mathbb{C}$ such that $T-\lambda id_H$ is not invertible. \\
Let $\mathcal{K}(H)$ be the compact operators on $H$, then the \textit{essential spectrum $\sigma_e(T)$ of T} is the spectrum of $T+\mathcal{K}(H)$ seen as element in the Calkin algebra $\mathcal{B}(H)/\mathcal{K}(H)$. Basic and important facts about both the spectrum and the essential spectrum are that they are bounded, closed, non-empty (if $\dim(H)=\infty)$ subsets of $\mathbb{C}$ and also that the essential spectrum is contained in the spectrum.\\
For the definition of the point and continuous spectrum, we have to take a closer look at when an operator is not invertible. It is obvious that this is the case if $\ker(T)\neq\{0\}$. Therefore we define the \textit{point spectrum $\sigma_p(T)$} to be all $\lambda\in\mathbb{C}$ such that $\ker(T-\lambda)\neq\{0\}$.\\
Now let $T^*$ be the adjoint of $T$, then the spectrum of $T$ can be decomposed into
  \begin{equation*}
    \sigma(T)=\sigma_p(T)\cup\overline{\sigma_p(T^*)}\cup\sigma_c(T),
  \end{equation*}
where the last set is defined as $\sigma_c(T)=\sigma(T)\setminus(\sigma_p(T)\cup\overline{\sigma_p(T^*)})$ and is called the \textit{continuous spectrum of $T$}. Using the equality $\ker(T^*-\bar\lambda id_H)=\textup{Im}(T-\lambda id_H)^\perp$, we see that a point $\lambda\in\mathbb{C}$ is in $\sigma_c(T)$ if and only if
  \begin{enumerate}[(i)]
    \item $\ker(T-\lambda id_H)=\ker(T^*-\bar\lambda id_H)=\{0\}$ and
    \item the image of $T$ is not closed.
  \end{enumerate}
Throughout the paper we will use (semi-)Fredholm operators. For completeness, we give the definition and the most important properties of (semi-)Fredholm operators for our purposes.\\ 
An operator $T\in\mathcal{B}(H)$ is called \textit{$($semi-$)$Fredholm} if 
  \begin{enumerate}[(i)]
    \item the image of $T$ is closed,
    \item $\dim(\ker(T))<\infty$ and (or) $\dim(\ker(T^*))<\infty$.
  \end{enumerate}
If $T$ is a semi-Fredholm operator, then $\textup{ind}(T)=\dim(\ker(T))-\dim(\ker(T^*))\in\mathbb{Z}\cup\{-\infty,\infty\}$ is called the \textit{index of T}. The following properties of semi-Fredholm operators are well known.

\begin{Lemma}
Let $c\in\mathbb{Z}\cup\{-\infty,\infty\}$. Then the set 
  \begin{equation*}
    \{T\in\mathcal{B}(H);\ T\ \textup{is semi-Fredholm with}\ \textup{ind}(T)=c\}
  \end{equation*}
is open with respect to the operator norm.
\label{Lemma1}
\end{Lemma}

\underline{Proof}: See for example \cite[Chapter 18 Corollary 2]{MV}.\\
\ \\
\begin{Lemma}
Let $T\in\mathcal{B}(H)$ and $\lambda\in\partial\sigma_e(T)$, then $T-\lambda id_H$ is not semi-Fredholm. In particular, if $\dim(H)=\infty$, there is always $\lambda\in\mathbb{C}$ such that $T-\lambda id_H$ is not semi-Fredholm.
\label{Lemma2}
\end{Lemma}

\underline{Proof}: The first claim is \cite[Chapter 19 Proposition 1]{MV}, the second claim follows from the fact that $\sigma_e(T)$ is compact and non-empty if $\dim(H)=\infty$.\\
\ \\
The last tool we need is the \textit{polar decomposition} of an operator and is defined as follows. Let $T\in\mathcal{B}(H)$ be an arbitrary operator, $P\in\mathcal{B}(H)$ a partial isometry and $A\in\mathcal{B}(H)$ a positive operator, then $P$ and $A$ are called the polar decomposition of $T$ if $\textup{Im}(A)^\perp=\ker(P)$ and $T=PA$. A proof of the existence and uniqueness is given in \cite[Problem 134]{HP}.\\ 
The next lemma and the following corollary are the starting point of the first two main theorems.

\begin{Lemma}
Let $T\in\mathcal{B}(H)$ be not semi-Fredholm and $\|T\|>\epsilon>0$. If $T=PA$ is the polar decomposition of $T$ and $\nu$ the spectral measure of $A$, then
  \begin{equation*}
    S=P\int_{\epsilon}^{\|T\|}t\ d\nu(t)\in\mathcal{B}(H)
  \end{equation*}
has the following properties:
  \begin{enumerate}[(i)]
    \item the image of $S$ is closed and more precisely $\|Sx\|\ge\epsilon\|x\|$ for all $x\in\ker(S)^\perp$,
    \item $\dim(\ker(S))=\dim(\textup{Im}(S)^\perp)=\infty$,
    \item $\|T-S\|\le\epsilon$.
  \end{enumerate}
\label{Lemma3}
\end{Lemma}

\underline{Proof}: Let $T, S, P, A, \nu$ be as above. Since $P$ is a partial isometry and $\nu([\epsilon,\|T\|])(H)\subset\nu(\{0\})(H)^\perp=\ker(A)^\perp=\overline{\textup{Im}(A)}=\ker(P)^\perp$, it applies that
  \begin{equation*}
    \left\|\int_{\epsilon}^{\|T\|}t\ d\nu(t)x\right\|=\|Sx\|
  \end{equation*}
for every $x\in \nu([\epsilon,\|T\|])(H)$. Therefore
  \begin{equation*}
    \|Sx\|^2=\left\|\int_{\epsilon}^{\|T\|}t\ d\nu(t)x\right\|^2=\left\langle \int_\epsilon^{\|T\|}t^2d\nu(t)x,x\right\rangle=\int_{\epsilon}^{\|T\|}t^2\ d\nu_{x,x}(t)\ge\epsilon^2\|x\|^2,
  \end{equation*}
where $\nu_{x,x}(\cdot)$ stands for the measure $\langle \nu(\cdot)x,x\rangle$, and together with $\ker(S)^\perp\subset\nu([\epsilon,\|T\|])(H)$, we have proven $(i)$.\\
For the sake of clarity, we denote by $[T]$ the equivalence class of $T+\mathcal{K}(H)$. Since $T$ is not semi-Fredholm, $[T]$ is not left-invertible by \cite[Theorem 14 and Theorem 15]{MV} and so $[T^*T]=[A]^2$ is not invertible. Thus $[A]$ is not invertible and in particular $A$ is not Fredholm. We can write $A$ as
  \begin{equation*}
    A=\int_{0}^{\|T\|}\textup{max}(\epsilon,t)d\nu(t)-\int_{0}^{\epsilon}(\epsilon-t)d\nu(t).
  \end{equation*}
The first operator $\int_{0}^{\|T\|}\textup{max}(\epsilon,t)d\nu(t)$ is obviously invertible
and if we assume that $\dim(\nu([0,\epsilon))(H))<\infty$, then the second operator $\int_{0}^{\epsilon}(\epsilon-t)d\nu(t)$ has finite rank. This contradicts the fact that $A$ is not Fredholm. Hence $\dim(\nu([0,\epsilon))(H))=\infty$. Together with
  \begin{equation*}
    \nu([0,\epsilon))(H)\subset \ker(S)\ \textup{and}\ P\nu([0,\epsilon))(H)\subset\textup{Im}(S)^\perp,
  \end{equation*}
where the second inclusion follows from $\nu([0,\epsilon))(H)\subset\ker(P)^\perp$, we obtain $\dim(\ker(S))=\dim(\textup{Im}(S)^\perp)=\infty$ and $(ii)$ is proven.\\
The last part $(iii)$ follows from 
  \begin{equation*}
    \|(T-S)x\|^2=\|\int t\chi_{[0,\epsilon)}(t)\ d\nu(t)x\|^2=\int t^2\chi_{[0,\epsilon)}(t)d\nu_{x,x}(t)\le\epsilon^2\|x\|^2
  \end{equation*}
for all $x\in H$.$\hfill\square$\\
\ \\

\begin{Corollary}
Let $T\in\mathcal{B}(H)$ and $\epsilon>0$. There are $\lambda\in\mathbb{C}$ and $S\in B_\epsilon(T-\lambda id_H)$ such that 
  \begin{enumerate}[(i)]
    \item the image of $S$ is closed,
    \item $\dim(\ker(S))=\dim(\textup{Im}(S)^\perp)=\infty$.
  \end{enumerate}
\label{Corollary1}
\end{Corollary}

\underline{Proof}: First apply \Cref{Lemma2} and then \Cref{Lemma3}.$\hfill\square$\\

\section{Non-empty point spectrum}
\label{Section 3}
Let $T$ be in $\mathcal{B}(H)$. It is well known that if $\dim(H)<\infty$, then every point in $\sigma(T)$ is an eigenvalue and so the point spectrum of $T$ is never empty. This is not true if $\dim(H)=\infty$. One of the best known examples is the bilateral shift on $\ell^2(\mathbb{Z})$. Now the question arises how $\glqq$big$\grqq$ the set of all operators with empty point spectrum is. Problem 8.4 in \cite{TE} asks if the set of operators with non-empty point spectrum is co-meager with respect to the operator norm. We will answer the question in the affirmative and even better, we will show that the complement of the set is \textit{nowhere dense}. 
\begin{Theorem}
The set 
  \begin{equation*}
    \{T\in\mathcal{B}(H); \ \sigma_p(T)=\emptyset\}
  \end{equation*}
is nowhere dense with respect to the operator norm.
\label{Theorem1}
\end{Theorem}
The proof follows immediately from the next lemma.

\begin{Lemma}
Let $T\in\mathcal{B}(H)$ and $\epsilon>0$. There exist $\delta>0$ and $\tilde T\in\mathcal{B}(H)$ such that
  \begin{enumerate}[(i)]
    \item $\|T-\tilde T\|<\epsilon$,
    \item $B_\delta(\tilde T)\subset \{A\in\mathcal{B}(H);\ \sigma_p(A)\neq\emptyset\}.$
  \end{enumerate}
\end{Lemma}

\underline{Proof}: We can assume that $\dim(H)=\infty$ or else $T=\tilde T$ fulfills $(i)$ and $(ii)$. Apply \Cref{Corollary1} to $T, \epsilon/2$ to obtain $S\in\mathcal{B}(H), \lambda\in\mathbb{C}$ and $\delta>0$ such that 
  \begin{enumerate}[(i)]
    \item the image of $S$ is closed,
    \item $\dim(\ker(S))=\dim(\ker(S^*))=\infty$,
    \item $\|S-(T-\lambda id_H)\|\le\epsilon/2$.
  \end{enumerate}
Since $\dim(\ker(S))=\dim(\ker(S^*))=\infty$, there is an isometric operator
  \begin{equation*}
    j:\ker(S^*)\to\ker(S)
  \end{equation*}
such that $\dim(\ker(S)\ominus\textup{Im}(j))<\infty$ and $j$ is not surjective. Let $I$ be the extension of $j$ to $H$ by setting $I(x)=0$ for $x\in\ker(S^*)^\perp=\textup{Im}(S)$. Then the operator
  \begin{equation*}
    \tilde{S}=S+\frac{\epsilon}{3}I^*
  \end{equation*}
is Fredholm with $\textup{index}(\tilde{S})=c<0$, which follows from $\textup{Im}(\tilde{S})=\textup{Im}(S)\oplus\textup{Im}(I^*)=H$ and $\ker(\tilde{S})=\ker(j^*)\neq\{0\}$. However, by \Cref{Lemma1} the set of Fredholm operators with index $c$ is open, therefore there is a $\delta>0$ with
  \begin{equation*}
    B_\delta(\tilde{S})\subset\{A\in\mathcal{B}(H);\ \textup{ind}(A)<0\}\subset\{A\in\mathcal{B}(H);\ \sigma_p(A)\neq\emptyset\}.
  \end{equation*}
The operator $\tilde{T}=\tilde{S}+\lambda id$ fulfills $(i)$, since
  \begin{equation*}
    \|T-\tilde{T}\|=\|T-\lambda id-\tilde{S}\|\le\|T-\lambda id-S\|+\|S-\tilde{S}\|<\epsilon,
  \end{equation*}
and also $(ii)$, since for $A\in B_\delta(\tilde{T})$, it holds that $A-\lambda id\in B_\delta(\tilde{S})$ and hence $\emptyset\neq\sigma_p(A-\lambda id)$. This immediately implies that $\emptyset\neq\sigma_p(A)$.
$\hfill\square$\\
\ \\
It is not new to study operators that are only determined by a spectral property with the help the polar decomposition. For example, Bouldin used the polar decomposition in \cite{Ri} to characterize the closure of the set of invertible operators.

\section{Non-empty continuous spectrum}
\label{Section 4}
As in the previous section, the case $\dim(H)<\infty$ is trivial, because then $\sigma_c(T)=\emptyset$. On the other hand, if $\dim(H)=\infty$, there are plenty of examples of operators with non-empty continuous spectrum, for example, every compact operator with trivial kernel. Now the question arose: Is the set of operators with empty continuous spectrum co-meager? (see \cite[Problem 8.4]{TE}) However, the following theorem answers the question in the negative. \\
\begin{Theorem}
The set 
  \begin{equation*}
    \{T\in\mathcal{B}(H);\sigma_c(T)\neq\emptyset\}
  \end{equation*}
is co-meager with respect to the operator norm.
\label{MainTheorem2}
\end{Theorem}

For the proof we need some technical lemmas and the Banach-Mazur game. We start with the definition of the later.\\
\ \\
Let $(X,\tau)$ be a topological space. Recall that a subset $U\subset X$ is called \textit{co-meager} if it contains a dense $G_\delta$ set, or equivalently there are sets $(U_n)_{n\in\mathbb{N}}$ in $X$ with
  \begin{enumerate}[(i)]
    \item $\textup{int}(\overline{X\setminus U_n})=\emptyset$,
    \item $U=\bigcap_{n\in\mathbb{N}}U_n$.
  \end{enumerate}
Now, assume two players, $I$ and $II$, are playing a game where they take turns playing sets of the topological space $(X,\tau)$
  \begin{equation*}
    \begin{split}
      &I:\ \ \ \ U_1\ \ \ \ \ \ \ \ U_3\ \ \dots\\
      &II:\ \ \ \ \ \ \ \ U_2\ \ \ \ \ \ \ \ U_4 \ \ \dots
    \end{split}
  \end{equation*}
with the condition that $\emptyset\neq U_{n+1}\subset U_n$ and $U_n\in\tau$ for all $n\in\mathbb{N}$. Player $II$ wins the game for a set $A\subset X$ if 
  \begin{equation*}
    \bigcap_{n\in\mathbb{N}}U_{2n}\subset A.
  \end{equation*}
Note that by choice of the $U_n$ we automatically have $\bigcap_{n\in\mathbb{N}}U_{2n-1}=\bigcap_{n\in\mathbb{N}}U_{2n}$.\\
The above game is called \textit{Banach-Mazur game} and the relation to co-meager sets is given by the next theorem. A proof can be found in \cite[Theorem 8.33]{AK}.

\begin{Theorem}
Let $(X,d)$ be a metric space and $\emptyset\neq A\subset X$. In the Banach-Mazur game, player $II$ has a winning strategy if and only if the set $A$ is co-meager.
\label{Banach-Mazur game}
\end{Theorem}

\begin{Remark}
If one has to show that Player $II$ has a winning strategy, it is enough to assume that Player $I$ only plays open balls $B_\epsilon(T)$. This follows from the condition that the sets played are open and Player $II$ has to find a winning strategy for all possible choices Player $I$ can make.
\label{Remark3.1}
\end{Remark}

Now we come to the technical lemmas needed to prove \Cref{MainTheorem2}.

\begin{Lemma}
Let $T\in\mathcal{B}(H)$ be not semi-Fredholm. Then for every $\epsilon>0$ there is $S\in B_\epsilon(T)$ such that 
  \begin{enumerate}[(i)]
    \item $\sigma_{e}(S)=\sigma_{e}(T)$,
    \item $\ker(S)=\ker(S^*)=\{0\}$.
  \end{enumerate}
\label{Lemma3.1}
\end{Lemma}

\underline{Proof}: The idea is to choose $S$ as a compact perturbation of $T$ since the essential spectrum is invariant under compact perturbations. Assume that there is a compact operator $K_1$ such that $\|K_1\|\le\epsilon/2$ and $\dim(\ker(T-K_1))=\dim(\ker(T^*-K_1^*))$. Then there is an injective, compact operator 
  \begin{equation*}
    K_2:\ker(T-K_1)\to\ker(T^*-K_1^*)
  \end{equation*}
with dense image and $\|K_2\|<\epsilon/2$. Extend $K_2$ to $H$ by $0$ on $\ker(T-K_1)^\perp$. Now the operator $S=T-K_1+K_2$ fulfills the conditions $(i)$, $(ii)$ and $S\in B_\epsilon(T)$, so it suffices to show that such an operator $K_1$ exists.\\
By \cite[Chapter 16 Theorem 18]{MV} there is a compact operator $K\in\mathcal{B}(H)$ such that $\dim(\ker(T-K))=\infty$ and $\|K\|<\epsilon/4$. We are done if $\dim(\ker(T^*-K^*))=\infty$. So assume that $\dim(\ker(T^*-K^*))<\infty$. Since $T^*-K^*$ is not semi-Fredholm, $\textup{Im}(T^*-K^*)$ is not closed. That is why we can again apply the same theorem to the operator
  \begin{equation*}
    H\to\ker(T-K)^\perp,x\mapsto (T^*-K^*)x.
  \end{equation*}
This gives another compact operator $\tilde K:H\to \ker(T-K)^\perp$ such that $\dim(\ker(T^*-K^*-\tilde K))=\infty$ and $\|\tilde K\|<\epsilon/4$. However, we can consider $\tilde K$ as an operator from $H$ to $H$. Finally, the operator $K_1=K+\tilde K^*\in\mathcal{B}(H)$ is compact as sum of two compact operators and fulfills $\|K_1\|<\epsilon/2$ as well as $\dim(\ker(T-K_1))=\dim(\ker(T^*-K_1^*))=\infty$. $\hfill\square$\\
\ \\

\begin{Lemma}
Let $T\in\mathcal{B}(H)$ be not semi-Fredholm with $\ker(T)=\{0\}$ and $T=PA$ the polar decomposition with spectral measure $v$. Then
  \begin{equation*}
    v([\delta,||T||])\to_{SOT}id_H\textup{\ for\ } \delta\to0.
  \end{equation*}
In particular, for every final dimensional subspace $F\subset H$, we have
  \begin{equation}
    ||(v([\delta,||T||])-id_H)|_F||\to0
\label{eq1}
  \end{equation}
for $\delta\to0$.
\label{Lemma3.2}
\end{Lemma}

\underline{Proof}:\\
It applies that $v((0,\|T\|])=id_H$, since $v((0,||T||])$ is the projection onto $\ker(A)^\perp$ and thus for each $x\in H$,
  \begin{equation*}
    \lim_{\delta\to0}\|(id_H-v([\delta,\|T\|))x\|^2=\lim_{\delta\to0}\|v([0,\delta))x\|^2=\lim_{\delta\to0}v_{x,x}([0,\delta))=v_{x,x}(\{0\})=0,
  \end{equation*}
where $v_{x,x}(\cdot)$ stands for the measure $\langle v(\cdot)x,x\rangle$. Hence $v([\delta,||T||])$ converges in the SOT to $id_H$ for $\delta\to0$. The additional note follows from a characterization of the SOT in \cite[ Problem 225]{HP}.$\hfill\square$\\
\ \\
A simple application of the triangle inequality shows that a sequence of norm bounded operators that converges pointwise on a dense subspace to an operator already converges pointwise on the entire space. This proves the next lemma.

\begin{Lemma}
Let $(v_n)_{n\in\mathbb{N}}$ be a sequence of spectral measures and $(r_n)_{n\in\mathbb{N}}$ in $[0,\infty)$ such that $v_n((0,r_n])=id_H$. Also assume that there is a increasing sequence of finite dimensional subspaces $F_n$ in  $H$, whose union is dense in $H$, and a null sequence $(\delta_n)_{n\in\mathbb{N}}$ of positive numbers such that
  \begin{equation*}
    ||(v_n([\delta_n,r_n])-id_H)|_{F_n}||<1/n.
  \end{equation*}
Then 
  \begin{equation*}
    v_n([\delta_n,r_n])x\to x
  \end{equation*}
for every $x\in H$.
\label{Lemma3.3}
\end{Lemma}
\ \\
\begin{Lemma}
Let $\epsilon>0, T\in\mathcal{B}(H), \lambda\in\partial\sigma_e(T)$. Then there is a $S\in B_{\epsilon}(T)$ and a $\delta>0$ such that
  \begin{equation*}
    B_\epsilon(\lambda)\cap\partial\sigma_e(A)\neq\emptyset
  \end{equation*}
for every $A\in B_\delta(S).$ 
\label{Lemma3.4}
\end{Lemma}
\ \\
\underline{Proof}:\\
Let $\epsilon>0$, $T\in\mathcal{B}(H)$ and $\lambda\in\partial\sigma_e(T)$. Since $\lambda$ is in the boundary of $\sigma_e(T)$, $T-\lambda id_H$ is not semi-Fredholm by \Cref{Lemma2} and additionally, there is a $\tilde\lambda\in B_{\epsilon}(\lambda)$ such that $T-\tilde\lambda id_H$ is Fredholm. Denote $\textup{ind}(T-\tilde\lambda id_H)=c$. By \Cref{Lemma1} the set of Fredholm operators with index $c$ is open, so there is a $\tilde\delta>0$ with
  \begin{equation}
    \textup{ind}(A-\tilde\lambda id_H)=c \ \forall\ A\in B_{\tilde\delta}(T).
\label{ind}
  \end{equation}
Without loss of generality, let $\tilde\delta<\epsilon$. Applying \Cref{Lemma3} to $T-\lambda id_H, \tilde\delta/2$ yields a $\tilde S\in\mathcal{B}(H)$ with the properties
  \begin{enumerate}[(i)]
    \item $\textup{Im}(\tilde S)$ is closed,
    \item $\textup{dim}(\ker(\tilde S))=\textup{dim}(\textup{Im}(\tilde S)^\perp)=\infty$,
    \item$ ||T-\lambda id_H-\tilde S||\le\tilde\delta/2$.
  \end{enumerate}
Take an isometric operator $j:\ker(\tilde S)\to\textup{Im}(\tilde S)^\perp$ such that $\dim(\ker(j))-\dim(\ker(j^*))=\tilde c\in\mathbb{Z}\setminus\{c\}$ and extend it to $H$ by defining $j(x)=0$ for $x\in\ker(\tilde S)^\perp$. Now the operator $S=\tilde S+\lambda id_H+\tilde\delta/3 j$ fulfills
  \begin{enumerate}[(a)]
    \item $||T-S||\le||T-\lambda id_H-\tilde S||+\tilde\delta/3||j||<\tilde\delta$,
    \item $S-\lambda id_H$ is Fredholm with $\textup{ind}(S-\lambda id_H)=\tilde c$,
    \item $S-\tilde\lambda id_H$ is semi-Fredholm with $\textup{ind}(S-\tilde\lambda id_H)=c$,
  \end{enumerate}
where part $(c)$ follows from $S\in B_{\tilde\delta}(T)$ and \Cref{ind}. By \Cref{Lemma1} there is a $\delta>0$ such that
  \begin{equation*}
    \begin{split}
      &\textup{ind}(A-\lambda id_H)=\tilde c \ \forall\ A\in B_{\delta}(S),\\
      &\textup{ind}(A-\tilde\lambda id_H)=c \ \forall\ A\in B_{\delta}(S).
    \end{split}
  \end{equation*}
It remains to show that $B_\epsilon(\lambda)\cap\partial\sigma_e(A)\neq\emptyset$ for every $A\in B_\delta(S).$ For this purpose, we define 
  \begin{equation*}
    t_0=\textup{min}\{t\in[0,1],\  A-\lambda id_H+t(\lambda-\tilde\lambda)id_H\ \textup{is not semi-Fredholm}\}
  \end{equation*}
and $\lambda_0=\lambda-t_0(\lambda-\tilde\lambda)$. Note that $t_0$ exists because $\textup{ind}( A-\lambda id_H)\neq\textup{ind}(A-\tilde \lambda id_h)$. By definition $A-\lambda_0 id_H$ is not semi-Fredholm, thus $\lambda_0\in\sigma_e(A)$ and since every neighborhood of $\lambda_0$ contains a point $\lambda_1$ such that $A-\lambda_1 id_H$ is Fredholm, we obtain $\lambda_0\in\partial\sigma_e(A).$ Finally, we still have to show that $\lambda_0$ is in $B_\epsilon(\lambda)$, but this follows from
  \begin{equation*}
    \pushQED{\qed} 
      |\lambda-\lambda_0|\le|\lambda-\tilde\lambda|<\epsilon. \qedhere
    \popQED
  \end{equation*}
\ \\
Now we have all the tolls and lemmas to prove \Cref{MainTheorem2}. For that fix an ONB $(e_n)_{n\in\mathbb{N}}$ of $H$ and define $F_n=\langle e_1,\dots,e_n\rangle$.\\
 \\
\underline{Proof of \Cref{MainTheorem2}}:\\
Recall that we want to use the Banach-Mazur game. By \Cref{Remark3.1} we can assume that player $I$ plays the set $B_{\epsilon_1}(T_1)\subset\mathcal{B}(H)$. Let $\lambda_1\in\partial\sigma_{e}(T_1)$, then by \Cref{Lemma3.1} we can assume that $\ker(T_1-\lambda_1id_H)=\ker(T_1^*-\bar\lambda_1id_H)=\{0\}.$\\
Let $T_1-\lambda_1 id_H=P_1A_1,\ T_1^*-\bar\lambda_1 id_H=Q_1B_1$ be the polar decompositions and $\nu_1, \mu_1$ be the spectral measures of $A_1, B_1$. By \Cref{Lemma3.2} there is a $1>\tilde{\delta_1}>0$ with
  \begin{equation*}
    \begin{split}
      &||(\nu_1([\tilde{\delta_1},||A_1||])-id_H)|_{F_1}||<1,\\
      &||(\mu_1([\tilde{\delta_1},||B_1||])-id_H)|_{F_1}||<1,
    \end{split}
  \end{equation*}
and by \Cref{Lemma3.4} there is a $\delta_1>0, S_1\in B_{\tilde{\delta}_1/4}(T_1)$ such that
  \begin{equation*}
    B_{\tilde{\delta}_1/4}(\lambda_1)\cap\partial\sigma_e(A)\neq\emptyset
  \end{equation*}
for every $A\in B_{\delta_1}(S_1)$. Without loss of generality, we can assume that $\delta_1<1$. Now Player II plays $B_{\delta_1}(S_1)$.\\
\ \\
Let $n>1.$ We will to construct the sets of Player II inductively. So assume that the sets $B_{\epsilon_i}(T_i), B_{\delta_i}(S_i), i=1,\dots,n-1$ have been played. Set $\lambda_0=0, \tilde\delta_0=\infty$. The induction hypothesis is :\\
There are $\lambda_i\in\partial\sigma_{e}(T_i)\cap B_{\tilde{\delta}_{i-1}/4}(\lambda_{i-1})$ with polar decompositions 
  \begin{equation*}
    T_i-\lambda_i id_H=P_iA_i,\ \ T_i^*-\bar\lambda_i id_H=Q_iB_i
  \end{equation*}
and spectral measures $\nu_i, \mu_i$ such that
  \begin{equation*}
    \begin{split}
      &||(\nu_i([\tilde{\delta_i},||A_i||])-id_H)|_{F_i}||<1/i, \\
      &||(\mu_i([\tilde{\delta_i},||B_i||])-id_H)|_{F_i}||<1/i
    \end{split}
  \end{equation*}
for $\tilde\delta_i>0$ with $1/i>\tilde\delta_i$ and $B_{\tilde{\delta}_i/4}(\lambda_i)\subset B_{\tilde{\delta}_{i-1}/4}(\lambda_{i-1})$. As well as $\delta_i<\min\{1/i,\tilde\delta_i/4\}$ and
  \begin{equation*}
    B_{\tilde{\delta}_i/4}(\lambda_i)\cap\partial\sigma_e(A)\neq\emptyset
  \end{equation*}
for every $A\in B_{\delta_i}(S_i)$.\\
 Assume Player I plays $B_{\epsilon_n}(T_n)$. By construction, there is a $\lambda_n\in\partial\sigma_{e}(T_n)\cap B_{\tilde{\delta}_{n-1}/4}(\lambda_{n-1})$ and by \Cref{Lemma3.1} we can assume that $\ker(T_n-\lambda_n)=\ker(T_n^*-\bar \lambda_n)=\{0\}$. Again, let $T_n-\lambda_n id_H=P_nA_n,\ T_n^*-\bar\lambda_n id_H=Q_nB_n$ be the polar decompositions and $v_n, \mu_n$ the spectral measures of $A_n, B_n$.
By \Cref{Lemma3.2} there is a $1/n>\tilde{\delta}_n>0$ such that
  \begin{equation*}
    \begin{split}
      &||(v_n([\tilde{\delta}_n,||A_n||])-id_H)|_{F_n}||<1/n,\\
      &||(\mu_n([\tilde{\delta}_n,||B_n||])-id_H)|_{F_n}||<1/n.\\
    \end{split}
  \end{equation*}
Without loss of generality, we can assume that $\tilde{\delta}_n$ is so small that $B_{\tilde{\delta}_n/4}(\lambda_n)\subset B_{\tilde{\delta}_{n-1}/4}(\lambda_{n-1})$.\\
By \Cref{Lemma3.4} there is a $\delta_n>0$ and a $S_n\in B_{\tilde{\delta}_n/4}(T_n)$ such that
  \begin{equation*}
    B_{\tilde{\delta}_n/4}(\lambda_n)\cap\partial\sigma_e(A)\neq\emptyset
  \end{equation*}
for every $A\in B_{\delta_n}(S_n)$. Without loss of generality, let $\delta_n<\min\{\tilde{\delta}_n/4,1/n\}$.\\
Now Player II plays $B_{\delta_n}(S_n)$.\\
\ \\
It remains to show that $T\in\bigcap_{n\in\mathbb{N}}B_{\delta_n}(S_n)=\bigcap_{n\in\mathbb{N}}B_{\epsilon_n}(T_n)$ fulfills $\sigma_c(T)\neq\emptyset$. \\
Let $\lambda_n, \tilde{\delta}_n, A_n, v_n$ be as above. Since $B_{\tilde{\delta}_i/4}(\lambda_i)\subset B_{\tilde{\delta}_{i-1}/4}(\lambda_{i-1})$ for every $i\in\mathbb{N}$, the sequence $(\lambda_n)_{n\in\mathbb{N}}$ converges to a $\lambda\in\mathbb{C}$ with $|\lambda-\lambda_n|\le\tilde\delta_n/4.$ On the other hand, $T_n$ converges to $T$ since $||T_n-T||<\epsilon_n\le\delta_{n-1}<1/(n-1).$ Thus $T-\lambda id_H$ is not semi-Fredholm because it is the limit of non semi-Fredholm operators.\\
It remains to show that $\ker(T-\lambda id_H)=\ker(T^*-\bar\lambda id_H)=\{0\}.$
Let $x\in\ker(T-\lambda id_H)$ and $x_n=v_n([\tilde{\delta}_n,||A_n||])x$. It applies that
  \begin{equation*}
    \begin{split}
      \tilde{\delta}_n(||x_n||-||x-x_n||)&\le||(T_n-\lambda_n id_H)x|||\\
      &=||(T_n-S_n)x+(S_n-T)x-(\lambda_n-\lambda)x||\\
      &\le\delta_n||x||+\delta_n||x||+\tilde{\delta}_n/4||x||\\
      &\le\tilde{\delta}_n/4||x||+\tilde{\delta}_n/4||x||+\tilde{\delta}_n/4||x||\\
      &=3/4\tilde{\delta}_n||x||.
    \end{split}
  \end{equation*}
Here the first inequality follows from $||(T_n-\lambda_n id_H)x_n||\ge\tilde\delta_n||x_n||$ and $||(T_n-\lambda_n id_H)(x-x_n)||\le\tilde\delta_n||x-x_n||$. By \Cref{Lemma3.3} $x_n$ converges to $x$ and so the above inequality shows that $x=0$.\\
One can do the same for $x\in\ker(T^*), x_n=\mu_n([\tilde\delta_n,\|A_n\|])x$ and obtains overall that $\ker(T-\lambda id_H)=\ker(T^*-\bar \lambda id_H)=\{0\}.$ Thus $\lambda\in\sigma_c(T).$
$\hfill\square$\\
\ \\
It is easy to see that the intersection of two co-meager sets is co-meager again. Therefore, \Cref{Theorem1} and \Cref{MainTheorem2} imply that having non-empty point and continuous spectrum is a typical property.

\section{Density of operators with empty continuous spectrum}
\label{Section 5}

Even with the knowledge that operators with empty continuous spectrum are meager, it is still not clear how big the closure of this set is. In this section we show that the closure of the set of operators with empty continuous spectrum is $\mathcal{B}(H)$:

\begin{Theorem}
Let $H$ be a separable Hilbert space. Then
  \begin{equation*}
    \overline{\{ T\in\mathcal{B}(H); \sigma_c(T)=\emptyset\}}=\mathcal{B}(H)
  \end{equation*}
with respect to the norm topology.
\label{MainTheorem3}
\end{Theorem}

The main idea of the proof is to show that every operator is in the closure of the similarity orbit of an operator with empty continuous spectrum. To achieve this, we will use the restricted Similarity Orbit Theorem from \cite{DH2}. \\
For an operator $T\in\mathcal{B}(H)$ we call 
  \begin{equation*}
    \textup{Sim}(T)=\{S^{-1}TS;\ S\in\mathcal{B}(H)\ \textup{invertible}\}
  \end{equation*}
the \textit{similarity orbit} of $T$. In \cite{DH}, the closure of $\textup{Sim}(T)$ was characterized in terms of properties of the spectrum of $T$. To formulate the restricted version of the Similarity orbit Theorem mentioned there, we need the following definitions.\\
The \textit{normal spectrum} $\sigma_0(T)$ are the isolated points in $\sigma(T)$ which are not in the essential spectrum $\sigma_e(T)$.\\
The points $\lambda\in\mathbb{C}$ such that $T-\lambda id_H$ is semi-Fredholm are denoted by $\rho_{sF}(T)$, and it immediately follows from \Cref{Lemma1} that this set is open.\\
Let $\lambda\in\sigma_e(T)$ be isolated and $\Phi$ be a faithful unital *-representation of the Calkin algebra $\mathcal{B}(H)/\mathcal{K}(H)$. Let $f(z)$ be the function that is equal to $z-\lambda$ on a neighborhood of $\lambda$ and $0$ on a neighborhood of $\sigma_e(T)\setminus\{\lambda\}$. We can apply the holomorphic functional calculus to the element $\Phi(T+\mathcal{K}(H))$ to obtain a quasinilpontent element $Q_\lambda$ . For $\lambda\in\mathbb{C}$, we now define
  \begin{equation*}
    k(\lambda;T)=
      \begin{cases}
        0 & \textup{if}\ \lambda\notin\sigma_e(T),\\
        n & \textup{if}\ \lambda\ \textup{is isolated in}\ \sigma_e(T)\ \textup{and} \ Q_\lambda\ \textup{is nilpotent of order} \ n,\\
\infty & \textup{otherwise}.
      \end{cases}
  \end{equation*}
It is not immediately clear that the above definition of $k(\lambda;T)$ is independent of the representation $\Phi$. However, this can be obtained by applying the function $f(z)$ to the element $T+\mathcal{K}(H)$. Then the element $f(T+\mathcal{K}(H))$ is nilpotent of order $n$ if and only if $k(\lambda;T)=n$. (see also \cite[Chapter 8.5]{DH1} and \cite[Chapter 9.1]{DH2}) The reason why we gave the definition that uses a representation of the Calkin algebra is that we use it in the proof of the next lemma.
\begin{Lemma}
Let $A\in\mathcal{B}(H)$ and $B\in\mathcal{B}(\tilde H)$ be operators on Hilbert spaces $H, \tilde H$. If $\lambda\in\sigma_e(A\oplus B)$ is isolated, then
  \begin{equation*}
    k(\lambda;A\oplus B)=\textup{max}(k(\lambda;A),k(\lambda;B)).
  \end{equation*}
\label{SumK}
\end{Lemma}
\underline{Proof:} Let $\lambda\in\sigma_e(A\oplus B)$ be isolated. For the sake of clarity, we denote the equivalence class of an operator $T$ in the Calkin algebra by $[T]$. Define the projection $P$ by
  \begin{equation*}
    P:H\oplus\tilde H\to H\oplus\tilde H, x\oplus y\mapsto x\oplus 0.
  \end{equation*}
Let $\Phi:\mathcal{B}(H\oplus \tilde H)/\mathcal{K}(H\oplus \tilde H)\to\mathcal{B}(K)$ be a faithful unital *-representation. Then $\Phi([P])$ is a projection, since $\Phi$ is a *-homomorphism. Denote this projection by $P_1$ and the Hilbert space $P_1(K)$ by $H_1$. Furthermore, we denote the projection $id_{K}-P_1$ by $P_2$ and $P_2(K)=H_2$. This yields a decomposition of $K$ into $H_1\oplus H_2$. Since $P_1$ commutes with $\Phi([A\oplus0])$ and $P_2$ with $\Phi([0\oplus B])$, we see that $H_1$ and $H_2$ reduce $\Phi([A\oplus B])$. Write $A_0$ for $\Phi([A\oplus B])|_{H_1}$ and $B_0$ for the restriction to $H_2$. Now
  \begin{equation*}
    \Phi_1:\mathcal{B}(H)/\mathcal{K}(H)\to \mathcal{B}(H_1), [T]\mapsto \Phi([T\oplus 0])|_{H_1}
  \end{equation*}
is well defined since $P_1$ commutes with $\Phi([T\oplus 0])$ for all $T\in\mathcal{B}(H)$ and also a faithful unital *-representation. Thus we obtain that $\Phi_1([A])=A_0$, $\sigma([A])=\sigma(A_0)$ and likewise for
  \begin{equation*}
    \Phi_2:\mathcal{B}(H)/\mathcal{K}(H)\to \mathcal{B}(H_2), [T]\mapsto \Phi([0\oplus T])|_{H_2},
  \end{equation*}
we obtain that $\Phi_2([B])=B_0$ and $\sigma([B])=\sigma(B_0)$. \\
Let $f(z)$ be a function that is the equal to $z-\lambda$ on a neighborhood of $\lambda$ and $0$ on a neighborhood of $\sigma_e(A\oplus B)\setminus\{\lambda\}$. We can apply $f$ to $A_0$ and $B_0$ since $\sigma_e(A\oplus B)=\sigma_e(A)\cup\sigma_e(B)$ to obtain
  \begin{equation*} 
    f(\Phi([A\oplus B]))=f(A_0\oplus B_0)=f(A_0)\oplus f(B_0).
  \end{equation*}
But $\Phi_1$ and $\Phi_2$ are faithful unital *-representations and since $k(\lambda;\cdot)$ is independent of the representation, we obtain that $f(A_0)$ is nilpotent of order $n$ if and only if $k(\lambda;A)=n$. Likewise $f(B_0)$ is nilpotent of order $n$ if and only if $k(\lambda;B)=n$. Thus $f(\Phi([A\oplus B]))$ is nilpotent of order $n$ if and only if $n=\textup{max}(k(\lambda;A),k(\lambda;B))$, and if $f(\Phi([A\oplus B]))$ is not nilpotent then either $k(\lambda;A)=\infty$ or $k(\lambda;B)=\infty$. We therefore have proven the lemma. $\hfill\square$\\
\ \\
Let $A$ be an operator in $\mathcal{B}(H)$. We say that $A$ has property $(S), (F)$ or $(A)$ with respect to $T$ if 
  \begin{itemize}
    \item[(S)]$ \sigma_0(A)\subset \sigma_0(T)$\ \textup{and each component of}\ $\mathbb{C}\setminus\rho_{sF}(A)$\ intersects\ $\sigma_e(T)$,\\
    \item[(F)] $\rho_{sF}(A)\subset\rho_{sF}(T),\ \textup{ind}(\lambda id_H-A)=\textup{ind}(\lambda id_H-T)$ and
      \begin{equation*}
        \min \textup{ind}(\lambda id_H-T)^k\le \min \textup{ind}(\lambda id_H-A)^k
      \end{equation*}
      for all $\lambda\in\rho_{sF}(T)$ and $k\ge1$,
 \\
    \item[(A)] $\dim(\ker(A-\lambda id_H))=\dim(\ker(T-\lambda id_H))$ for all $\lambda\in \sigma_0(A)$.
  \end{itemize}
Here $\min \textup{ind} (x)$ stands for
  \begin{equation*}
    \min\{\dim(\ker(x)),\dim(\ker(x^*))\}.
  \end{equation*}
Now we are ready to formulate the restricted Similarity Orbit Theorem \cite[Corollary 1.6]{DH}.

\begin{Theorem}[restricted Similarity Orbit Theorem]
Let $T\in\mathcal{B}(H)$ and $k(\lambda;T)=\infty$ for every isolated point $\lambda\in\sigma_e(T)$. Then
  \begin{equation*}
    \overline{\textup{Sim}(T)}=\{X\in\mathcal{B}(H); X\ \textup{satisfies property}\ (S), (F)\ \textup{and}\ (A)\ \textup{with respect to}\ T\}.
  \end{equation*}
\label{SimOrb}
\end{Theorem}

We need the following special operator to get rid of the points $\lambda\in\sigma_e(A)$ with $k(\lambda;A)\in\mathbb{N}$ .

\begin{Example}
Let $T_0$ be a quasinilpotent operator in $\mathcal{B}(H)$ such that $k(0;T_0)=\infty$ and
 $a=(\lambda_n)_{n\in\mathbb{N}}$ be a bounded sequence of isolated points in $\mathbb{C}$. Then the operator
  \begin{equation*}
    T_a:\bigoplus_{n\in\mathbb{N}} H\to\bigoplus_{n\in\mathbb{N}} H, (x_n)_n\mapsto ((\lambda_n id_H-T_0)(x_n))_{n}
  \end{equation*} 
is in $\mathcal{B}(\bigoplus_{n\in\mathbb{N}}H)$ since the sequence $a$ is bounded. Furthermore if $\lambda\notin\overline{\{\lambda_n,\ n\in\mathbb{N}\}}$, then, by continuity of the inverse map, $\sup_{n\in\mathbb{N}}\|((\lambda-\lambda_n)id_H-T_0)^{-1}\|<\infty$ and hence $\lambda\notin\sigma(T_a)$. We obtain that $\overline{\{\lambda_n;\ n\in\mathbb{N}\}}=\sigma(T_a)=\sigma_e(T_a)$ and each $\lambda_n$ is a isolated point in $\sigma_e(T_a)$. Fix $m\in\mathbb{N}$. We can write $T_a$ as a direct sum of $\lambda_m id_H-T_0$ and 
  \begin{equation*}
    \bigoplus_{n\in\mathbb{N}\setminus\{m\}}H\to\bigoplus_{n\in\mathbb{N}\setminus\{m\}}H, (x_n)_n\mapsto ((\lambda_n id_H-T_0)(x_n))_{n}.
  \end{equation*}
By assumption, we have that $k(\lambda_m;\lambda_m id_H-T_0)=\infty$, but $\lambda_m$ is isolated in $\sigma_e(T_a)$ and so we conclude that $k(\lambda_m;T_a)=\infty$ by \Cref{SumK} . Since $m$ was arbitrary, we see that $k(\lambda_n;T_a)=\infty$ for all $n\in\mathbb{N}$.
\label{ExQN}
\end{Example}

\begin{Lemma}
Let $\tilde H$ be a Hilbert space and $A\in\mathcal{B}(H), B\in\mathcal{B}(\tilde H)$. Assume that $\sigma(B)=\sigma_p(B)=\mathbb{C}\setminus\rho_{sF}(A)$. Then
  \begin{enumerate}[(i)]
    \item $\sigma(A\oplus B)=\sigma(A)$,
    \item $\sigma_c(A\oplus B)=\emptyset$,
    \item $\sigma_e(A\oplus B)=\sigma_e(A)$,
    \item $\sigma_0(A\oplus B)=\sigma_0(A)$ and $\dim(\ker(A\oplus B-\lambda id_{H\oplus\tilde H}))=\dim(\ker(A-\lambda id_H))$ for all $\lambda\in\sigma_0(A)$,
    \item $\rho_{sF}(A\oplus B)=\rho_{sF}(A)$,
    \item $\textup{ind}(A\oplus B-\lambda id_{H\oplus\tilde H})=\textup{ind}(A-\lambda id_{H})$ for all $\lambda\in\rho_{sF}(A\oplus B)$ ,
    \item $\min \textup{ind}(A\oplus B-\lambda id_{H\oplus \tilde H})^k=\min \textup{ind}(A-\lambda id_H)^k$ for all $\lambda\in\rho_{sF}(A\oplus B)$ and $k\ge1$.
  \end{enumerate}
If $U$ is a linear, invertible operator from $H\oplus \tilde H\to H$, then $A$ has property $(S), (A)$ and $(F)$ with respect to $U(A\oplus B)U^{-1}$.
\label{LemmaUnionSpectra}
\end{Lemma}

\underline{Proof}: $(i)$ follows from $\sigma(A\oplus B)=\sigma(A)\cup\sigma(B)$ and $\sigma(B)=\mathbb{C}\setminus\rho_{sF}(A)\subset\sigma(A)$.\\
Let $\lambda\in\sigma_c(A\oplus B)$. $\lambda$ is in $\sigma(A)$ but not in $\rho_{sF}(A)$, so we get that $\lambda\in\sigma(B)$. However, $\sigma(B)=\sigma_p(B)\subset\sigma_p(A\oplus B)$. This is a contradiction to $\lambda\in\sigma_c(A\oplus B)$ and we obtain $(ii)$.\\
For $(iii)$ observe that $\sigma_e(A\oplus B)=\sigma_e(A)\cup\sigma_e(B)$ and $\sigma_e(B)\subset\sigma(B)\subset\sigma_e(A)$.\\
From $(i)$ and $(iii)$ we get that $\sigma_0(A\oplus B)=\sigma_0(A)$. In addition, we have that 
  \begin{equation}
    \dim(\ker(A\oplus B-\lambda id_{H\oplus\tilde H}))=\dim(\ker(A-\lambda id_H))+\dim(\ker(B-\lambda id_{\tilde H}))
\label{eqdimker}
  \end{equation}
for all $\lambda\in\mathbb{C}$ and together with $\dim(\ker(B-\lambda id_{\tilde H}))=0$ for all $\lambda\in\sigma_0(A)\subset\mathbb{C}\setminus\sigma(B)$, we obtain $(iv)$.\\
$(v)$ and $(vi)$ essentially follow from $\sigma(B)\subset\mathbb{C}\setminus\rho_{sF}(A)$ and $\textup{ind}(B-\lambda id_{\tilde H})=0$ for all $\lambda\in\rho_{sF}(A)$. To be more precise, $\textup{Im}(A\oplus B-\lambda id_{H\oplus\tilde H})$ is closed if and only if the image of $A-\lambda id_H$ and $B-\lambda id_{\tilde H}$ is closed, and together with \Cref{eqdimker}, we see that if $A-\lambda id_H$ is semi-Fredholm and $B-\lambda id_{\tilde H}$ is Fredholm, then $A\oplus B-\lambda id_{H\oplus\tilde H}$ is semi-Fredholm and fulfills
  \begin{equation*}
    \textup{ind}(A\oplus B-\lambda id_{H\oplus\tilde H})=\textup{ind}(A-\lambda id_H)+\textup{ind}(B-\lambda id_{\tilde H})
  \end{equation*}
for all $\lambda\in\rho_{sF}(A\oplus B)$. Since $B-\lambda id_{\tilde H}$ is invertible for all $\lambda\in\rho_{sF}(A)$, it applies that $\textup{ind}(B-\lambda id_{\tilde H})=0$ and $\rho_{sF}(A)\subset\rho_{sF}(A\oplus B)$. However, if $A-\lambda id_H$ is not semi-Fredholm, then $A\oplus B-\lambda id_{H\oplus \tilde H}$ is not semi-Fredholm. This shows the reversed inclusion $\rho_{sF}(A\oplus B)\subset\rho_{sF}(A)$.\\
Part $(vii)$ follows from 
  \begin{equation*}
    \begin{split}
      \ker(A\oplus B-\lambda id_{H\oplus\tilde H})=\ker(A-\lambda id_{H})\oplus\ker(B-\lambda id_{\tilde H}), \\
      (A\oplus B-\lambda id_{H\oplus\tilde H})^k=(A-\lambda id_H)^k\oplus(B-\lambda id_{\tilde H})^k
    \end{split}
  \end{equation*}
 and $\ker((B-\lambda id_{\tilde H})^*)^k=\ker(B-\lambda id_{\tilde H})^k=\{0\}$ for all $\lambda\in\rho_{sF}(A)=\rho_{sF}(A\oplus B)$.\\
The additional remark follows from $(i), \dots, (viii)$ and the observation that each of the sets $\sigma(\cdot), \sigma_e(\cdot), \sigma_c(\cdot), \rho_{sF}(\cdot), \sigma_0(\cdot)$ and each number $\dim(\ker(\cdot)), \textup{ind}(\cdot), \min\textup{ind}(\cdot)$ is invariant under similarity.$\hfill\square$\\
\ \\
As a last result, before we give the proof of \Cref{MainTheorem3}, we need that for every compact set $K\subset\mathbb{C}$, there is an operator $T$ on a separable Hilbert space such that 
  \begin{equation*}
    \sigma(T)=\sigma_p(T)=K.
  \end{equation*}
Fortunately, this has already been shown by Dixmier and Foia\c{s} \cite{ JD}. In fact, the question of the topological properties of the spectrum and the point spectrum has a long history. Independently of Dixmier and Foia\c{s}, Nikol'skaya showed in \cite{Nik} that the point spectrum of a bounded operator acting on a separable reflexive Banach space is $F_\sigma$ and that any bounded $F_\sigma$ subset of $\mathbb{C}$ coincides with the point spectrum of a certain bounded operator acting on a separable Hilbert space. A few years later, Kaufmann showed in \cite{Ka1, Ka2} that a necessary and sufficient condition for a subset of $\mathbb{C}$ to be the point spectrum of a bounded operator in some separable complex Banach space is that it be analytic in the sense of Suslin and bounded. All these results were extended by Smolyanov and Shkarin. For further results we recommend \cite{SaS} and \cite{Nie}.\\
\\
\underline{Proof of \Cref{MainTheorem3}}:\\
We already know that the theorem is true for $\dim(H)<\infty$, we therefore assume for the rest of the proof that $\dim(H)=\infty$. Let $A\in\mathcal{B}(H)$ and $\{\lambda_n; n\in\mathbb{N}\}$ be the set of isolated points in $\sigma_e(A)$. If the set is finite we simply repeat a $\lambda$ infinitely often. Let $T_a$ be the operator from \Cref{ExQN} with respect to the sequence $(\lambda_n)_{n\in\mathbb{N}}$. Recall that $\sigma(T_a)=\overline{\{\lambda_n;\ n\in\mathbb{N}\}}$ and $k(\lambda_n;T_a)=\infty$ for all $n\in\mathbb{N}$. By \cite{JD} and since the set $\mathbb{C}\setminus\rho_{sF}(A)$ is compact and non-empty, there is an operator $T_p\in\mathcal{B}(H)$ such that 
  \begin{equation*}
    \sigma(T_p)=\sigma_p(T_p)=\mathbb{C}\setminus\rho_{sF}(A).
  \end{equation*}
Now define the operators $B=T_p\oplus T_a\in\mathcal{B}(H\oplus\bigoplus_{n\in\mathbb{N}}H)$ and $\tilde T=A\oplus B$. Fix a linear, invertible operator $U:H\oplus H\oplus\bigoplus_{n\in\mathbb{N}} H\to H$ and define $T=U\tilde TU^{-1}$.\\
The set 
  \begin{equation*}
    \{\lambda\in\rho_{sF}(A);\ \textup{ind}(A-\lambda id_H)=\infty\ \textup{or}\ \textup{ind}(A-\lambda id_H)=-\infty\}
  \end{equation*}
is open by \Cref{Lemma1} and equal to $\sigma_e(A)\cap\rho_{sF}(A)$. Hence no isolated point in $\sigma_e(A)$ can belong to $\rho_{sF}(A)$ and therefore $\{\lambda_n,\ n\in\mathbb{N}\}$ is contained in $\mathbb{C}\setminus\rho_{sF}(A)$. It also applies that $\sigma(B)=\sigma(T_a)\cup\sigma(T_p)=\mathbb{C}\setminus\rho_{sF}(A)$ and $\sigma_p(B)=\sigma_p(T_a)\cup\sigma_p(T_p)=\sigma(B)$. Now we can use \Cref{LemmaUnionSpectra} to obtain that $A$ has property $(A), (S)$ and $(F)$ with respect to $T$. In addition, \Cref{SumK} and the invariance of $k(\lambda;\cdot)$ under similarity, imply that $k(\lambda_n;T)=k(\lambda_n;T_a)=\infty$ for all $n\in\mathbb{N}$. Thus the requirements for the restricted Similarity Orbit Theorem \ref{SimOrb} are fulfilled and we obtain that $A\in\overline{\textup{Sim}(T)}$. But since $T$ has empty continuous spectrum, we see that 
  \begin{equation*}
    \pushQED{\qed} 
A\in\overline{\textup{Sim}(T)}\subset\overline{\{ X\in\mathcal{B}(H);\ \sigma_c(X)=\emptyset\}}.\qedhere
    \popQED
  \end{equation*}

The idea of constructing an operator with certain spectral properties by direct sums is well known and appears, for example, in the proof of \cite[Theorem 2]{SaS}.\\
We end the paper with a theorem that characterises the closure of all operators with the property that the point spectrum is equal to the spectrum. The proof is only a slight modification of the proof of \Cref{MainTheorem3}.

\begin{Theorem}
It applies that
  \begin{equation*}
    \overline{\{T\in\mathcal{B}(H);\ \sigma_p(T)=\sigma(T)\}}=\left\{T\in\mathcal{B}(H);
      \begin{split} 
        \dim(\ker(T-\lambda id_H))\neq0\ \\
        \ \textup{for all}\ \lambda\in\rho_{sF}(T)\cap\sigma(T)
      \end{split}
    \right\}.
  \end{equation*}
\end{Theorem}
\underline{Proof}: Let $T\in\mathcal{B}(H)$ and assume that there is a $\lambda\in\rho_{sF}(T)\cap\sigma(T)$ such that $\dim(\ker(T-\lambda id_H))=0$. Since the operator $T-\lambda id_H$ is semi-Fredholm and has trivial kernel, it is bounded below by some constant $c>0$. By \Cref{Lemma1} there is a $\epsilon<c/2$ such that $\textup{ind}(S)=\textup{ind}(T-\lambda id_H)\neq0$ for all $S\in B_\epsilon(T-\lambda id_H)$. Now each operator $S\in B_\epsilon(T)$ fulfills that $S-\lambda id_H$ is bounded below and $\lambda\in\sigma(S)$. Hence 
  \begin{equation*}
    B_\epsilon(T)\subset\mathcal{B}(H)\setminus\{X\in\mathcal{B}(H);\ \sigma_p(X)=\sigma(X)\}
  \end{equation*}
and in particular, we have that $T\notin\overline{\{X\in\mathcal{B}(H);\ \sigma_p(X)=\sigma(X)\}}$.\\
On the other hand, if $A\in\mathcal{B}(H)$ and $\dim(\ker(A-\lambda id_H))\neq0$ for all $\lambda\in\rho_{sF}(A)\cap\sigma(A)$, then we only have to check that the operator $T$ from the proof of \Cref{MainTheorem3} lies in $\{X\in\mathcal{B}(H);\ \sigma_p(X)=\sigma(X)\}$. However, by construction of $T$, we have that $\mathbb{C}\setminus\rho_{sF}(T)\subset\sigma_p(T),\ \sigma(T)=\sigma(A),\ \rho_{sF}(T)=\rho_{sF}(A)$ and by the assumptions on the operator $A$, we have that $\rho_{sF}(A)\cap\sigma(A)\subset\sigma_p(A)$. Since $\sigma_p(A)\subset\sigma_p(T)$, we obtain $\sigma_p(T)=\sigma(T)$ and the proof is complete. $\hfill\square$

\bibliography{BibTypicalSpectra} 
\bibliographystyle{plain}
Fachrichtung Mathematik, Universit\"at des Saarlandes, 66123 Saarbr\"ucken, Germany\\
\textit{Email address:} scherer@math.uni-sb.de

\end{document}